\documentclass[11pt]{article}
\usepackage{amsmath,amsfonts,amssymb,hyperref,setspace,cite}
\newtheorem{thm}{Theorem}[section]
\newtheorem{lem}[thm]{Lemma}
\newtheorem{cor}[thm]{Corollary}
\newtheorem{exa}[thm]{Example}
\newtheorem{defn}[thm]{Definition}

\numberwithin{equation}{section}
\title{ Contractions Which Obtained  by Equivalent Metric of Cone Metric}
\author{\\{\normalsize{\sc
 Mehdi Asadi\footnote{Corresponding author.\quad Fax:+98-241-4220030,\quad Email: masadi.azu@gmail.com}}}\\
   \\{\it Department of Mathematics, Zanjan Branch,}
   \\{\it Islamic Azad University,  Zanjan, Iran}
  \\{\it { masadi@iauz.ac.ir}}   }
        \date{\empty}
\begin{document}
\maketitle

\begin{abstract}
In this paper we try to collect certain contractions which can be obtained by equivalent metric of cone metric.\\

\end{abstract}

\section{Introduction and Preliminary}
 Long-Guang and Xian in \cite{Xian} generalized the concept of a metric space, replacing the set of real numbers by an ordered Banach space and obtained some fixed point theorems for mapping satisfy different contractive conditions.\\
Recently Wei-Shih Du in \cite{Du1} has proved that the Banach contraction principle in general metric spaces and in TVS-cone metric spaces are equivalent, and in \cite{Du2} has obtained new type fixed point theorems for nonlinear multivalued maps in metric spaces and the generalizations of Mizoguchi-Takahashi's fixed point theorem and Berinde-Berinde's fixed point theorem. But in this paper according to metric which introduced by Feng and Mao in \cite{Feng}, Asadi and Vaezpour in \cite{AVS-3}, we obtain the equivalent contractive conditions which satisfies as their contractions in cone metrics.\\

Let $E$ be a real Banach space. A nonempty convex closed subset
$P\subset E$ is called a cone in $E$ if it satisfies:
\begin{enumerate}
  \item[$(i)$] {$P$ is closed, nonempty and $P\neq \{0\}$,}
  \item[$(ii)$] {$a,b\in \mathbb{R},$ $a,b\geq 0$ and $x,y \in P$ imply that $ax+by \in P,$}
  \item[$(iii)$] {$x \in P$ and $-x \in P$ imply that $x = 0.$}
\end{enumerate}
The space $E$ can be partially ordered by the cone $P\subset E$;
that is, $x \le y$ if and only if  $y-x \in P$. Also we write
$x\ll y$ if $y-x \in P^o$, where $P^o$ denotes the interior of $P$.\\
A cone $P$ is called normal if there exists a constant $k>0$ such
that $0\le x \le y$ implies $\|x\| \le k\|y\|$.\\
In the sequel we always suppose that $E$ is a real Banach
space, $P$ is a cone in $E$ with nonempty interior i.e. $P^o\neq \emptyset$ and $\leq$ is the partial ordering with
respect to $P$.
\begin{defn}
(\cite{Xian}) Let $X$ be a nonempty set.
Assume that the mapping $D:X\times X\rightarrow E$ satisfies
\begin{enumerate}
  \item[(i)] {$0\leq D(x,y)$ for all $x,y \in X$ and $D(x,y)=0 $ iff $x=y$}
  \item[(ii)] {$D(x,y)=D(y,x)$ for all $x,y \in X$}
  \item[(iii)] {$D(x,y)\leq D(x,z)+D(z,y)$ for all $x,y,z \in X$}.
\end{enumerate} Then $D$ is called a cone metric on $X$, and
$(X,D)$ is called a cone metric space.
\end{defn}
\begin{defn}
Metric $d$ is equivalent of cone metric $D$, if generated topology of $d$ and $D$ be equal and further the equivalent metric satisfies the same contractive conditions as the cone metric.
\end{defn}
In other words, convergence one of them implies that other ones, i.e.
$$x_n \overset{d}{\longrightarrow} x\iff x_n \overset{D}{\longrightarrow} x.$$
\begin{thm}(\cite{Feng})
For every cone metric $D:X\times X\rightarrow E$ there exists a metric $d:X\times X\rightarrow \mathbb{R}^+$ which is equivalent to $D$ on $X$.
\end{thm}
Indeed, the metric $d$ that has been defined in \cite{Feng,AVS-3} is
$d(x,y)=\inf\{\|u\|: D(x,y)\leq u\}.$
     Also, remember that for all $\{x_n\}\subseteq X$ and $x\in X$, $x_n\rightarrow x$ in $(X,d)$ if and only if $x_n\rightarrow x$ in $(X,D)$ (\cite{Feng,AVS-3}).\\
\newline
Throughout this paper we shall show that the equivalent metric satisfies the same contractive conditions as cone metric. So most of the fixed point theorems which have been proved are the straightforward results from the metric case. For more details see \cite{HSh,AR,IR1,IR2,R.Sh,Wardowski,deim,Rh,RR,KRR,JKRR,JRRR,R,KRV,Rh2}.
 \section{Main Results}
\begin{lem}\label{b}
Let $D,D^*:X\times X\rightarrow \Bbb E$ be cone metrics, $d,d^*:X\times X\rightarrow \mathbb{R}^+$ their equivalent metrics respectively and $T:X\rightarrow X$ a self map. If  $D(Tx,Ty)\leq D^*(x,y)$, then $d(Tx,Ty)\leq d^*(x,y).$
\end{lem}
{\bf Proof.} By the definition of $d^*$,
$$\forall \varepsilon>0 \quad\exists v \quad \|v\|<d^*(x,y)+\varepsilon,\quad D^*(x,y)\leq v.
$$
Therefore if $D(Tx,Ty)\leq D^*(x,y)\leq v,$
then we have $$d(Tx,Ty)\leq \| v\|\leq d^*(x,y)+\varepsilon.$$
Since $\varepsilon>0$ was arbitrary so $d(Tx,Ty)\leq d^*(x,y).$ $\Box$\newline
\begin{exa}
Let $E:=\Bbb R^+$, $P:=\Bbb R^+$ and $D:X\times X\rightarrow E$ be a cone metric, $d:X\times X\rightarrow \mathbb{R}^+$ be its equivalent metric. Also let  $T:X\rightarrow X$ be a self map and $\varphi:\Bbb R^+\rightarrow \Bbb R^+$ be defined by $\varphi (x)=\frac{x}{1+x}$. If  $D^*:=\varphi(D)$, then
it is easy to see that $D^*(x,y)=\varphi (D(x,y))=\frac{D(x,y)}{1+D(x,y)}$ is a cone metric and its equivalent metric is $d^*=\varphi(d)$, and if, $D(Tx,Ty)\leq \varphi (D(x,y))=\frac{D(x,y)}{1+D(x,y)},$
then by Lemma \ref{b}, $d(Tx,Ty)\leq \varphi (d(x,y))=\frac{d(x,y)}{1+d(x,y)}.$
We can see that $x_n\to x$ in $(X,d)$ if and only if $x_x\to x$ in $(X,D)$.
\end{exa}
\begin{defn}
A self map $\varphi$ on a normed space $X$ is bounded if
$$\|\varphi\|:=\sup_{0\neq x\in X}\frac{\|\varphi(x)\|}{\|x\|}<\infty.$$
\end{defn}
\begin{thm}
Let $D:X\times X\rightarrow  E$ be a cone metric, $d:X\times X\rightarrow \mathbb{R}^+$ its equivalent metric, $T:X\rightarrow X$ a self map and $\varphi:P\rightarrow P$ a bounded map, then there exists
$\psi:\mathbb{R}^+\rightarrow \mathbb{R}^+$
such that $D(Tx,Ty)\leq \varphi(D(x,y))$ for every $x,y\in X$ implies $d(Tx,Ty)\leq \psi(\|D(x,y)\|)$ for all $x,y\in X$. Moreover if $\psi$ is decreasing map or $\varphi$ is linear and increasing map then,
$d(Tx,Ty)\leq \psi(d(x,y))$ for all $x,y\in X$.
\end{thm}
{\bf Proof.} Put $\psi(t):=\sup_{0\neq x\in P}\left\|\varphi\left(\frac{t}{\|x\|}x\right )\right\|$ for all $t\in \Bbb R^+$ and note that $\psi(t)\leq t\|\varphi\|$ for all $t\in\Bbb R^+$.
So $\|\varphi(x)\|\leq\psi(\|x\|)$  for all $x\in P$.
Therefore if $D(Tx,Ty)\leq \varphi(D(x,y)),$
then we have $d(Tx,Ty)\leq \| \varphi(D(x,y))\|\leq \psi(\|D(x,y)\|).$
\newline
By the definition of $d$ we have $d(x,y)\leq \|D(x,y)\|$. Now if $\psi$ is a decreasing map,then
$$d(Tx,Ty)\leq  \psi(\|D(x,y)\|)\leq \psi(d(x,y)).$$  If $\varphi$ is a
linear increasing map, then $\psi(t)=t\|\varphi\|.$ The definition of $d$ implies that
$$\forall \varepsilon>0 \quad\exists v \quad \|v\|<d(x,y)+\varepsilon,\quad D(x,y)\leq v.
$$
Therefore if $D(Tx,Ty)\leq \varphi(D(x,y))\leq \varphi(v),$
then we have
$$d(Tx,Ty)\leq \| \varphi(v)\|\leq \psi(\|v\|)\leq \psi(d(x,y))+\psi(\varepsilon).$$
Since $\varepsilon>0$ was arbitrary and $\psi(\varepsilon)\rightarrow 0$ as $\varepsilon\rightarrow 0$,  so $d(Tx,Ty)\leq  \psi(d(x,y)).$ $\Box$\newline
\par
In the following summary of our results are listed.

\begin{cor}\label{2.5}
Let $D$ be a cone metric, $d$ its equivalent metric, $T :X\rightarrow X$ a map, $\lambda\in[0,\frac{1}{2})$ and $\alpha,\beta\in [0,1)$. For $x, y \in X$,
\begin{enumerate}
\item[i.] $D(Tx, Ty)\leq \alpha D(x,y) \Rightarrow d(Tx, Ty)\leq \alpha d(x, y).$
  \item[ii.] $D(Tx, Ty)\leq \lambda (D(Tx, x) + D(Ty, y)) \Rightarrow d(Tx, Ty)\leq \lambda (d(Tx, x) + d(Ty, y)).$
  \item[iii.] $D(Tx, Ty)\leq \lambda (D(Tx, y) + D(Ty, x)) \Rightarrow d(Tx, Ty)\leq \lambda (d(Tx, y) + d(Ty, x)).$
  \item[iv.] $D(Tx, Ty)\leq \alpha D(x, Ty) + \beta D(Tx, y) \Rightarrow d(Tx, Ty)\leq \alpha d(x, Ty) +\beta d(Tx, y).$
      \item[v.] $D(Tx, Ty)\leq \alpha D(x, Tx) + \beta D(y, Ty) \Rightarrow d(Tx, Ty)\leq \alpha d(x, Tx) +\beta d(y, Ty).$
\end{enumerate}
   \end{cor}
\begin{cor}
Let $D$ be a cone metric, $d$ its equivalent metric, $T :X\rightarrow X$ a map and $\alpha,\beta\in [0,1)$. For $x, y \in X$,
\begin{enumerate}
      \item[a.] there exists $u\in \{D(x, y); D(x, Tx); D(y, Ty); \frac{1}{2}[
D(x, Ty)]+ D(y, Tx)]\}$ such that $D(Tx, Ty)\leq \alpha u$ where $\alpha \in(0, 1),$
 then $$d(Tx, Ty)\leq \alpha\max\{d(x, y); d(x, Tx); d(y, Ty); \frac{1}{2}[d(x, Ty)]+ d(y, Tx)]\}.$$
  \item[b.] there exists $u\in \{D(x, y); D(x, Tx); D(y, Ty); \frac{1}{2}
D(x, Ty);\frac{1}{2}D(y, Tx)\}$ such that $D(Tx, Ty)\leq \beta u$ where $\beta \in(0, 1),$ then $$d(Tx, Ty)\leq \beta\max\{d(x, y); d(x, Tx); d(y, Ty); \frac{1}{2}d(x, Ty);\frac{1}{2}d(y, Tx)\}.$$
  \item[c.] there exists $u\in \{D(x, y); \frac{1}{2}[D(x, Tx)+ D(y, Ty)]; \frac{1}{2}[D(x, Ty)+ D(y, Tx)]\}$ such that $D(Tx, Ty)\leq \beta u$ where $\beta \in(0, 1),$ then $$d(Tx, Ty)\leq \beta\max\{d(x, y); \frac{1}{2}[d(x, Tx)+ d(y, Ty)]; \frac{1}{2}[d(x, Ty)+ d(y, Tx)]\}.$$
  \end{enumerate}
     \end{cor}
{\bf Proof.} To prove $(a)$, if $u\in \{D(x, y); D(x, Tx); D(y, Ty)\}$, then by Corollary \ref{2.5}, $(i)$; and if $u= \frac{1}{2}(D(x, Ty)+ D(y, Tx))$  by Corollary \ref{2.5}, $(iv)$; with $\alpha=\beta=\frac{1}{2}$ we obtain desire results.\newline
$(b)$ and $(c)$ are clear, by Corollary \ref{2.5}, $(i)$ and $(iv),(v)$ respectively.  $\Box$
  \begin{cor}
Let $D$ be a cone metric, $d$ its equivalent metric, $T :X\rightarrow X$ a map. For $x, y \in X$,
\begin{enumerate}
  \item[a.] if $$D(Tx, Ty)\leq a_1D(x, y)+a_2 D(x, Tx)+a_3 D(y, Ty)+a_4D(x, Ty)+a_5D(y, Tx),$$ then $$d(Tx, Ty)\leq a_1d(x, y)+a_2 d(x, Tx)+a_3 d(y, Ty)+a_4d(x, Ty)+a_5d(y, Tx)$$ where $\sum_{i=1}^5a_i< 1.$
    \item[b.] if there exists $$u\in \{D(x, y); D(x, Tx); D(y, Ty); D(x, Ty);D(y, Tx)\}$$ such that  $D(Tx, Ty)\leq \lambda u,$ then $$d(Tx, Ty)\leq \lambda\max\{d(x, y); d(x, Tx); d(y, Ty); d(x, Ty);d(y, Tx)\}$$ where $\lambda\in[0,\frac{1}{2})$.
  \item[c.] if $$D(Tx, Ty)\leq a_1D(x, y)+a_2 D(x, Tx)+a_3 D(y, Ty)+a_4[D(x, Ty)+D(y, Tx)],$$ then $$d(Tx, Ty)\leq a_1d(x, y)+a_2 d(x, Tx)+a_3 d(y, Ty)+a_4[d(x, Ty)+d(y, Tx)]$$ where $a+1+a_2+a_3+2a_4< 1.$
     \end{enumerate}
          \end{cor}
{\bf Proof.} To prove $(a)$, for convenience, put
$$D^T:=D(Tx,Ty),D_1:=D(x,y),$$
$$D_2:=D(x,Tx),D_3:=D(y,Ty),D_4:=D(x,Ty),D_5:=D(y,Tx)$$
and similarly
$$d^T:=d(Tx,Ty),d_1:=d(x,y),$$
$$d_2:=d(x,Tx),d_3:=d(y,Ty),d_4:=d(x,Ty),d_5:=d(y,Tx).$$
So $D^T\leq \sum_{i=1}^5a_iD_i$. Now by definition of $d$
$$\forall i~(1\leq i\leq 5)~\forall \varepsilon>0~\exists v_i\quad\text{s.t.} \quad \|v_i\|<d_i+\varepsilon$$
and $D_i\leq v_i$. Therefore
$$D^T\leq \sum_{i=1}^5a_iD_i\leq \sum_{i=1}^5a_iv_i,$$
thus
$$d^T\leq \|\sum_{i=1}^5a_iv_i\|\leq \sum_{i=1}^5a_i\|v_i\|< \sum_{i=1}^5a_id_i+(\sum_{i=1}^5a_i)\varepsilon,$$
since $\varepsilon>0$ is arbitrary so we have
$$d^T\leq \sum_{i=1}^5a_id_i.$$
To prove  $(b)$ and $(c)$ we use the Corollary \ref{2.5}. $\Box$
\begin{cor}
Let $D,D^*$ be cone metrics, $d,d^*$ their equivalent metrics, $T :X\rightarrow X$ a map. There exist $m,n\in\mathbb{N}$ and $ k\in[0,1)$ such that $$D(T^mx,T^ny)\leq kD(z,t)$$ for all $x,y\in X$, $z\neq t$ and $z,t\in \{x,y,T^px,T^qy\}$ where $1\leq p\leq m$ and $1\leq q\leq n,$ then $$d(T^mx,T^ny)\leq k d(z,t).$$
\end{cor}
{\bf Proof.} By definition of $d$ we have
$$\forall z,t\in \{x,y,T^px,T^qy\}~\forall \varepsilon>0~\exists v\quad\text{s.t.} \quad \|v\|<d(z,t)+\varepsilon$$
where $z\neq t$ and $D(z,t)\leq v$.  So
$$D(T^mx,T^ny)\leq kD(z,t)\leq kv,$$ therefore
$$d(T^mx,T^ny)\leq \|kv\|< kd(z,t)+k\varepsilon$$
since $\varepsilon>0$ is arbitrary so we have
$$d(T^mx,T^ny)\leq k d(z,t). \Box$$
{\bf Acknowledgements}\newline
The author is indebted to referee for carefully reading the paper and for making useful suggestions. This paper has been supported by the Zanjan Branch, Islamic Azad University, Zanjan, Iran. The author would like to thanks this support.



\begin{thebibliography}{99}
\bibitem{Xian} Long-Guang, Z. Xian, Cone metric spaces and fixed point theorems of contractive mapping, J. Math. Anal. Appl. 322(2007),  1468--1476.
\bibitem{Du1} Wei-Shih Du, A Note on Cone Metric Fixed Point Theory and its equivalence, Nonlinear Analysis, 72(2010), 2259-2261.
\bibitem{Du2} Wei-Shih Du, New Cone Fixed Point Theorems for Nonlinear Multivalued Maps with their Applications, Applied Mathematics Letters, article in press.(2010).
\bibitem{Feng} Y. Feng and W. Mao, Equivalence of Cone Metric spaces
and Metric Spaces, Fixed Point Theory, {\bf 11}(2)(2010), 259-264.
\bibitem{AVS-3} M. Asadi, S. M. Vaezpour,  Scalarization of contractions in the cone metric spaces, The First Workshop in Fixed Point Theory and Applications, Amirkabir University of Technology and IPM, Iran, June 09-10, 2010.
\bibitem{HSh} R. H. Haghi, S. Rezapour, Fixed points of multifunctions on regular cone metric  spaces, Expositiones Mathematicae, 28(1)(2010), 71-77.
\bibitem{AR} M. Abbas, B. E. Rohades, Fixed point and periodic point results in cone metric spaces, Appl. Math. Lett. 22(2009), 511-515.
\bibitem{IR1} D. Illi\'{c}, V. Rako\v{c}evi\'{c}, Common fixed point for maps on cone metric spaces, J.  Math. Anal. and Appl. 341(2008), 876-882.
\bibitem{IR2} D. Illi\'{c}, V. Rako\v{c}evi\'{c}, Quasi contraction on a cone metric spaces, Appl. Math.  Lett. 22(2009), 728-731.
\bibitem{R.Sh} S. Rezapour, R. Hamlbarani, Some notes on the paper cone metric spaces and  fixed point theorems of contractive mappings, J. Math. Anal. Appl. 345(2008), 719-724.
\bibitem{Wardowski} D. Wardowski, Endpoints and fixed points of set-valued contractions in cone  metric spaces, Nonliner Analysis: Theory, Methods and Applications 71(2009), 512-516.
\bibitem{deim} K. Deimling, Nonlinear Functional Analysis, Springer-Verlage, Germany, 1985.
\bibitem{Rh} B. E. Rhoades, A comparison of various definition of contractive mappings, Trans. Amer. Math. Soc. 266(1977), 256-290.
\bibitem{RR} S. Radenovi\'{c}, B.E. Rhoades, Fixed point theorem for two non-self mappings in cone metric spaces, Computers and Mathematics with Applications 57(2009), 1701-1707.
\bibitem{KRR} Z. Kadelburg, S. Radenovi\'{c} and B. Rosi\'{c}, Strict Contractive Conditions and Common Fixed Point Theorems in Cone Metric Spaces, Fixed Point Theory and Applications, 2009(2009), Article ID 173838, 14 pages.
\bibitem{JKRR} S. Jankovi\'{c}, Z. Kadelburg, S. Radenovi\'{c} and B. E. Rhoades, Assad-Kirk-Type Fixed Point Theorems for a Pair of Nonself Mappings on Cone Metric Spaces, Fixed Point Theory and Applications,  2009(2009), Article ID 761086, 16 pages.
\bibitem{JRRR} G. Jungck, S. Radenovi\'{c}, S. Radojevi\'{c} and V. Rako\v{c}evi\'{c}, Common Fixed Point Theorems for Weakly Compatible Pairs on Cone Metric Spaces, Fixed Point Theory and Applications, 2009 (2009), Article ID 643840, 13 pages.
\bibitem{R} S. Radenovi\'{c},  Common fixed points under contractive conditions in cone metric spaces, Computers and Mathematics with Applications 58(2009),  1273-1278.
\bibitem{KRV} Z. Kadelburg, S. Radenovi\'{c} and V. Rako\v{c}evi\'{c}, Remarks on "Quasi-contraction on a cone metric space", Appl. Math. Lett. 22(2009), 1674-1679.
\bibitem{Rh2} B. E. Rhoades, A fixed point theorem for some non-self-mappings, Mathematica Japonica 23(1978), 457-459.
\end{thebibliography}
\end{document}